\documentclass{article}

\usepackage{delarray,verbatim,enumerate,a4wide}
\usepackage{amsmath,amsthm,amstext,amsbsy,amssymb,amsfonts,amscd}
\usepackage[all]{xy}
\usepackage[french]{babel}
\usepackage[T1]{fontenc}
\usepackage[utf8]{inputenc}

\usepackage{scalerel}

\usepackage{upgreek,bbm}

\usepackage[small]{titlesec}

\usepackage{color}
\definecolor{vert}{rgb}{0.1,0.4,0.2}
\usepackage[colorlinks=true,linkcolor=blue,citecolor=vert]{hyperref}

\usepackage{calligra}
\DeclareFontShape{T1}{calligra}{m}{n}{<->s*[0.95]callig15}{}
\DeclareMathAlphabet{\mathscr}{T1}{calligra}{m}{n}

\newtheorem{Th}{Théorème}[]
\newtheorem{Lem}[Th]{Lemme}
\newtheorem{Prop}[Th]{Proposition}
\newtheorem{Cor}[Th]{Corollaire}

\newtheorem{Def} [Th]{Définition}

\newtheorem*{Th*}{Théorème}
\newtheorem*{Cor*}{Corollaire}
\newtheorem*{Def*}{Définition}

\def\Preuve{\noindent {\it Preuve.~}}

\def\Remarque{\smallskip\noindent {\it Remarque.~}}

\def\Exemple{\smallskip\noindent {\it Exemple.~}}

		\def\QQ{\mathbb Q}	
\def\NN{\mathbb N}	\def\ZZ{\mathbb Z}		
\def\F2{\mathbb{F}_2}	\def\Z2{\mathbb{Z}_2}		
\def\Zl{{\mathbb{Z}_\ell}} 	\def\Ql{\mathbb{Q}_\ell}	\def\Tl{\mathbb{T}_\ell}	

 		\def\P{\mathcal  P}		\def\U{\mathcal  U}	\def\F{\mathcal  F}	
  		\def\C{\mathcal  C}		\def\R{\mathcal  R}	\def\X{\mathcal  X}
 	  	\def\Cl{\mathcal  C\!\ell}	
\def\E{\mathcal  E}		\def\T{\mathcal  T}			\def\D{\mathcal  D}

		\def\mmu{\boldsymbol{\upmu}}

		\def\p{{\mathfrak p}}				\def\a{{\mathfrak a}}		\def\m{{\mathfrak m}}
\def\r{{\mathfrak R}}		\def\l{{\mathfrak l}}

	\def\deg{\operatorname{deg}}		\def\Ind{\operatorname{Ind}}
\def\Gal{\operatorname{Gal}}			\def\Rad{\operatorname{Rad}}
		\def\Hom{\operatorname{Hom}}		

\newcommand\scale[2]{\vstretch{#1}{\hstretch{#1}{#2}}}

\newcommand\si[1]{\scale{.7}{#1}}	
\newcommand\ph{{\phantom{*}}}
\newcommand\ab{{\scale{.8}{\rm ab}}}		
	\newcommand\reg{{\scale{.8}{\rm r\acute eg}}}	
		
\def\%{{\scale{.8}{\infty}}}

\makeatletter
\newcommand*\wt[2][0.2ex]{%
        \begingroup
        \mathchoice{\wt@helper{#1}{#2}{\displaystyle}{\textfont}}
                   {\wt@helper{#1}{#2}{\textstyle}{\textfont}}
                   {\wt@helper{#1}{#2}{\scriptstyle}{\scriptfont}}
                   {\wt@helper{#1}{#2}{\scriptscriptstyle}{\scriptscriptfont}}%
        \endgroup
        #2%
}
\newcommand*\wt@helper[4]{%
        \def\currentfont{\the#41}%
        \def\currentskewchar{\char\the\skewchar\currentfont}%
        \setbox\tw@\hbox{\currentfont$#2$\currentskewchar}%
        \dimen@ii\wd\tw@
        \setbox\tw@\hbox{\currentfont$#2${}\currentskewchar}%
        \advance\dimen@ii-\wd\tw@
        \rlap{\raisebox{-#1}{$\m@th#3\kern\dimen@ii\widetilde{\phantom{#2}}$}}%
}
\makeatother

\begin{document}

\title{\Large\bf Sur les Modules d'Iwasawa $S$-ramifiés $T$-décomposés}
\author{ Jean-François {\sc Jaulent} }
\date{}

\maketitle

{\footnotesize
\noindent{\bf Résumé.}
Nous corrigeons les formules fautives contenues dans un article précédent et explicitons le module de défaut pour les invariants $\lambda$ d'Iwasawa attachés aux pro-$\ell$-extensions abéliennes $S$-ramifiées $T$-décomposées sur la $\Zl$-extension cyclotomique d'un corps de nombres. Les formules obtenues recoupent et prolongent les résultats de Itoh, Mizusawa et Ozaki sur les modules d'Iwasawa modérément ramifiés.
}\smallskip

{\footnotesize
\noindent{\bf Abstract.}
We correct the faulty formulas given in a previous article and we compute the defect group for the Iwasawa $\lambda$ invariants attached to the $S$-ramified $T$-decomposed abelian pro-$\ell$-extensions over the $\Zl$-cyclotomic extension of a number field. As a consequence, we extend the results of Itoh, Mizusawa and Ozaki on tamely ramified Iwasawa modules for the cyclotomic $\Zl$-extension of abelian fields.}\smallskip


\section*{Introduction}
\addcontentsline{toc}{section}{Introduction}

Supposons donnés un corps de nombres $K$, un premier $\ell$ et une $\Zl$-extension $K_\%$ de $K$.\smallskip

Le résultat emblématique de la théorie d'Iwasawa (cf. e.g. \cite {Ser,Was}) affirme que les ordres respectifs $\ell ^{x(n)}$ des $\ell$-groupes de classes d'idéaux $\,\Cl(K_n)$  attachés aux étages finis $K_n$ de la tour $K_\infty/K$, de degrés respectifs $[K_n:K]=\ell^n$  sont donnés pour $n$ assez grand par la formule:\smallskip

\centerline{$x(n) \ =\ \mu \ell ^n \ + \ \lambda n \ + \ \nu$,}\smallskip

\noindent où $\nu$ est un entier relatif, et où $\lambda$ et $\mu$ sont des entiers naturels déterminés par la pseudo-décomposition de la limite projective (pour les applications normes) $\,\C(K_\%) =  \varprojlim  \, \Cl(K_n)$, regardée comme module de torsion sur l'algèbre d'Iwasawa $\Lambda = \Zl [[\gamma -1]]$ construite sur un générateur topologique $\gamma$ du groupe procyclique $\Gamma = \Gal(K_\infty/K)$.

\begin{Def*} Étant données $(a_n)_{n\in\NN}$ et $(b_n)_{n\in\NN}$ deux suites d'entiers relatifs, convenons d'écrire:\smallskip

\centerline{$a_n \simeq b_n$, lorsque la différence $a_n-b_n$ est bornée; $a_n \approx b_n$, lorsqu'elle est ultimement constante.}
\end{Def*}

Cela posé, l'égalité précédente, focalisée sur les seuls paramètres structurels $\lambda$ et $\mu$ s'écrit:\smallskip

\centerline{$x(n) \ \approx \ \mu \ell ^n \ + \ \lambda n$,}\smallskip

\noindent et vaut identiquement si l'on remplace les $\ell$-groupes $\,\Cl(K_n)$ par leurs quotients respectifs d'exposant $\ell ^n$ (ou $\ell^{n+k}$, pour $k$ fixé), comme expliqué dans \cite{J18}.\medskip

Soient maintenant $\bar S$ et $\bar T$ deux ensembles finis disjoints de places de $K$; et soit  $\Cl^{\si{\bar S}}_{\si{\bar T}}(K_n)$ le pro-$\ell$-groupe des $\bar T$-classes $\bar S$-infinitésimales de $K_n$. Ce pro-$\ell$-groupe correspond, par la théorie $\ell$-adique du corps de classes (cf. \cite{J31}), à la pro-$\ell$-extension abélienne maximale de $K_n$ qui est non-ramifiée en dehors des places divisant celles de $\bar S$ et totalement décomposée aux places au-dessus de celles de $\bar T$; et c'est en particulier un $\Zl$-module de type fini. Son quotient d'exposant $\ell^n$, disons ${}^{\ell^n}\!\Cl_{\si{\bar T}}^{\si{\bar S}}(K_n)$, est ainsi un $\ell$-groupe; et on s'attend à ce que la $\ell$-valuation ${x_{\si{\bar T}}^{\si{\bar S}}(n)}$ de son ordre s'exprime asymptotiquement de façon simple à partir des invariants structurels du module d'Iwasawa:\smallskip

\centerline{$\C_{\si{\bar T}}^{\si{\bar S}}(K_\%)= \varprojlim \,\Cl_{\si{\bar T}}^{\si{\bar S}}(K_n)$.}\smallskip

C'est le programme initié dans \cite{J18}, puis développé dans \cite{J40}. La formule attendue\smallskip

\centerline{$x_{\si{\bar T}}^{\si{\bar S}}(n) \ \simeq \ \rho_{\si{\bar T}}^{\si{\bar S}} n\ell^n \ + \ \mu_{\si{\bar T}}^{\si{\bar S}} \ell ^n \ + \ \lambda_{\si{\bar T}}^{\si{\bar S}} n$,}\smallskip

\noindent qui fait intervenir la quantité $\rho_{\si{\bar T}}^{\si{\bar S}}=\dim_\Lambda\,\C_{\si{\bar T}}^{\si{\bar S}}(K_\%)$ (i.e.  la dimension sur le corps des fractions $\Phi$ de $\Lambda$ du $\Phi$-espace $\Phi \otimes_\Lambda \C_{\si{\bar T}}^{\si{\bar S}}(K_\%)$) ainsi que la $\ell$-valuation $\mu_{\si{\bar T}}^{\si{\bar S}}$ et le degré $\lambda_{\si{\bar T}}^{\si{\bar S}}$ du polynôme caractéristique du sous-module de $\Lambda$-torsion $\T_{\si{\,\bar T}}^{\si{\bar S}}(K_\%)$ de $\, \C_{\si{\bar T}}^{\si{\bar S}}$ est cependant en défaut dans certains cas, comme repéré par Salle \cite{Sal}. Plus précisément, les paramètres $ \rho^{\si{\bar S}}_{\si{\bar T}}$et $ \mu_{\si{\bar T}}^{\si{\bar S}}$ coïncident bien avec les invariants structurels, mais ce n'est pas toujours le cas pour les paramètres $ \lambda_{\si{\bar T}}^{\si{\bar S}}$.\smallskip

Le résultat principal de \cite{J52}, corrigeant \cite{J40}, est le suivant:

\begin{Th*}[Jaulent--Maire--Perbet]\label{JMB}
Il existe un entier relatif $\kappa_{\si{\bar T}}^{\si{\bar S}}$ tel que l'on ait le paramétrage:\smallskip

\centerline{$x_{\si{\bar T}}^{\si{\bar S}}(n) \ \simeq \ \rho_{\si{\bar T}}^{\si{\bar S}} n\ell^n \ + \ \mu_{\si{\bar T}}^{\si{\bar S}} \ell ^n \ + \ (\lambda_{\si{\bar T}}^{\si{\bar S}} -\kappa_{\si{\bar T}}^{\si{\bar S}}) n$,}
\end{Th*}

Le but de la présente note est de rectifier les formules défectueuses de \cite{J31}, énoncées en termes de caractères, et de donner en particulier une formulation correcte des identités du miroir de Gras pour les caractères $\lambda_{\si{\bar T}}^{\si{\bar S}}$ attachés à la $\Zl$-extension cyclotomique d'un corps de nombres. Puis, dans un second temps, d'étudier le caractère de défaut $\kappa_{\si{\bar T}}^{\si{\bar S}}$ correspondant et d'en donner une interprétation arithmétique simple. Enfin, au moins sous certaines hypothèses, de le déterminer explicitement.\smallskip

Pour ne pas alourdir cette étude, nous nous  limitons au cas  au cas $\ell\ne 2$, techniquement plus facile, mais le cas $\ell=2$ relève essentiellement des mêmes méthodes. Nous supposons que $K$ est une extension abélienne contenant $\mmu_\ell^\ph$ d'un sous-corps totalement réel $F$ , de groupe de Galois $\Delta$ d'ordre étranger à $\ell$.
L'algèbre $\ell$-adique $\Zl [\Delta ]$ est alors une algèbre semi-locale, produit direct des extensions non ramifiées $Z_\varphi=\Zl[\Delta]e_\varphi$ de $\Zl$; et les idempotents primitifs $e_\varphi$ sont donnés à partir des caractères $\ell$-adiques irréductibles $\varphi $ de $\Delta $ par les formules classiques:\smallskip

\centerline{$e_\varphi\, =\, \frac{1}{d}\sum_{\tau\in \Delta}\varphi(\tau^{-1})\tau$.}\smallskip

De façon toute semblable, l'algèbre de groupe $\Lambda[\Delta]$ s'écrit canoniquement comme produit:\smallskip

\centerline{$\Lambda[\Delta] = \oplus_\varphi \,\Lambda[\Delta] e_\varphi = \oplus_\varphi\, \Lambda_\varphi\ $.}\smallskip

\noindent Et, pour chaque caractère irréductible $\varphi$ du groupe $\Delta$, la $\varphi$-composante $\Lambda_\varphi = \Lambda[\Delta] e_\varphi$ associée à l'idempotent $e_\varphi$ s'identifie à l'algèbre des séries formelles $Z_\varphi [[\gamma -1]]$ en l'indéterminée $\gamma -1$.\medskip

Plus généralement, par action des idempotents primitifs $e_\varphi$ tout $\Lambda[\Delta]$-module noethérien $X$ se décompose naturellement comme somme directe de ses $\varphi$-composantes $X_\varphi = X^{e_\varphi}$, chaque $X_\varphi$ étant pseudo-isomorphe, comme $\Lambda_\varphi$-module noethérien à un unique $\Lambda_\varphi$ module élémentaire:\smallskip

\centerline{$X_\varphi \sim \Lambda_\varphi^{\rho_\varphi} \oplus \big( \oplus_{i=0}^{s_\varphi}\Lambda_\varphi /f_{\varphi,i}\Lambda_\varphi \big) \oplus \big( \oplus_{j=0}^{t_\varphi}\Lambda_\varphi / \ell^{m_{\varphi,i}}\Lambda_\varphi \big)$.}\smallskip

\noindent Notant alors $P_\varphi = \prod_{j=0}^{t_\varphi}\ell^{m_{\varphi,j}}\prod_{i=0}^{s_\varphi}f_{\varphi,i} \ \in \ZZ_\varphi[\gamma-1]$ le polynôme caractéristique du sous-module de $\Lambda_\varphi$-torsion de $X_\varphi$, on obtient ainsi les trois invariants structurels $\rho_\varphi=\dim_{\Lambda_\varphi}X_\varphi$, $\mu_\varphi=\nu_\ell(P_\varphi)$ et $\lambda_\varphi=\deg P_\varphi$,
qu'il est commode de coder globalement en introduisant les {\em caractères structurels}:\smallskip

\centerline{$\rho = \sum_\varphi \rho_\varphi \ \varphi \ , \qquad
\mu = \sum_\varphi \mu_\varphi \ \varphi \ , \qquad
\lambda = \sum_\varphi \lambda_\varphi \ \varphi $.}\medskip

Appliquant cette construction au $\Lambda[\Delta]$-module noethérien $\,\C_{\si{\bar T}}^{\si{\bar S}}(K_\%)$ lorsque $K_\%=KF_\%$ provient d'une $\Zl$-extension arbitraire de $F$, on définit ainsi les trois caractères structurels $\rho^{\si{\bar S}}_{\si{\bar T}}$, $\mu^{\si{\bar S}}_{\si{\bar T}}$ et $\lambda^{\si{\bar S}}_{\si{\bar T}}$. Puis, en transposant {\em mutatis mutandis} le résultat de \cite{J52} rappelé plus haut aux $\varphi$-composantes des groupes $\,\Cl^{\si{\bar S}}_{\si{\bar T}}(K_n)$, on obtient immédiatement:

\begin{Th} [Théorème des paramètres]\label{TPa}
Soient $\ell$ un nombre premier impair et $K/F$ une extension abélienne de corps de nombres, de degré $d$ étranger à $\ell$, puis $F_\%=\bigcup F_n$ une $\Zl$-extension arbitraire de $F$ et $K_\%=\bigcup K_n$, avec $[K_n:K]=[F_n:F]=\ell^n$, de sorte qu'on a: $\Delta=\Gal(K_\%/F_\%)\simeq\Gal(K/F)$. Soit enfin $\gamma$ un générateur topologique du groupe procyclique $\Gamma=\Gal(K_\%/K)\simeq\Gal(F_\%/F)$ et $\Lambda=\Zl[[\gamma-1]]$ l'algèbre d'Iwasawa associée.\par

Deux ensembles finis disjoints $\bar S$ et $\bar T$ de places de $F$ étant donnés, notons $\rho^{\si{\bar S}}_{\si{\bar T}}$, $\mu^{\si{\bar S}}_{\si{\bar T}}$ et $\lambda^{\si{\bar S}}_{\si{\bar T}}$ les caractères structurels pour sa structure de $\Lambda[\Delta]$-module attachés à la limite projective (pour la norme) $\,\C^{\si{\bar S}}_{\si{\bar T}}(K_\%)=\varprojlim\,\Cl^{\si{\bar S}}_{\si{\bar T}}(K_n)$ des (pro)-$\ell$-groupes de $\bar T$-classes $\bar S$-infinitésimales des corps $K_n$.\smallskip

Il existe alors un caractère $\ell$-adique virtuel $\kappa^{\si{\bar S}}_{\si{\bar T}}$ de $\Delta$ tel que la $\ell$-valuation $x^{\si{\bar S}}_{\si{\bar T}}(n)_\varphi$ de la $\varphi$-composante du quotient d'exposant $\ell^n$ de $\,\Cl^{\si{\bar S}}_{\si{\bar T}}(K_n)$ soit donnée asymptotiquement par l'identité:\smallskip

\centerline{$x^{\si{\bar S}}_{\si{\bar T}}(n)_\varphi \ \approx \ <\rho^{\si{\bar S}}_{\si{\bar T}}, \varphi>n\ell^n \ + \ <\mu^{\si{\bar S}}_{\si{\bar T}}, \varphi>\ell^n \ + \ <\lambda^{\si{\bar S}}_{\si{\bar T}}-\kappa^{\si{\bar S}}_{\si{\bar T}}, \varphi>n $.}
\end{Th}

Nous nous proposons dans ce qui suit de préciser le {\em caractère de défaut} $\kappa^{\si{\bar S}}_{\si{\bar T}}$ lorsque le corps $F$ est totalement réel,  $F_\%$ sa $\Zl$-extension cyclotomique et $K/F$ à conjugaison complexe, i.e. pour $K$ extension quadratique totalement imaginaire d'un sur-corps $\bar K$ de $F$ totalement réel.\medskip

\Remarque Le produit scalaire $<\varphi,\varphi>$ est le degré $\deg\varphi=[Z_\varphi:\Zl]$ du caractère $\ell$-adique irréductible $\varphi$. Il vaut $1$ si et seulement si $\varphi$ est absolument irréductible; par exemple pour $d \mid (\ell-1)$.

\newpage
\section{Énoncé du théorème principal}

Précisons d'abord quelques conséquences immédiates de \cite{J52} dans le contexte général où l'on ne suppose ni que $F$ est totalement réel ni que $F_\%$ est sa $\Zl$-extension cyclotomique:

\begin{Prop}
Plaçons-nous dans la situation du Théorème \ref{TPa}. Alors:
\begin{itemize}
\item[(i)] Dans le cas spécial où l'extension $F_\infty/F$ est elle-même $\bar S$-ramifiée et $\bar T$-décomposée, on a $\kappa^{\si{\bar S}}_{\si{\bar T}}=-1$ (opposé du caractère unité) ainsi que l'estimation asymptotique stricte:\smallskip

\centerline{$x^{\si{\bar S}}_{\si{\bar T}}(n)_\varphi \ \approx \ <\rho^{\si{\bar S}}_{\si{\bar T}}, \varphi>n\ell^n \ + \ <\mu^{\si{\bar S}}_{\si{\bar T}}, \varphi>\ell^n \ + \ <\lambda^{\si{\bar S}}_{\si{\bar T}}+1, \varphi>n $.}\smallskip

\item[(ii)] Dans tous les autres cas, on a: $0 \leqslant \kappa^{\si{\bar S}}_{\si{\bar T}} \leqslant \ell^e\,\rho^{\si{\bar S}}_{\si{\bar T}} $, où $e$ est le plus petit entier $n$ tel que les places ramifiées dans $F_\%/F$ le sont totalement dans $F_\%/F_n$ et les places de $T$ finiment décomposées dans $F_\%/F_n$ ne se décomposent pas dans 
$F_\%/F_n$.
En particulier, pour $\varphi$ fixé $<\kappa^{\si{\bar S}}_{\si{\bar T}},\varphi>$ est nul dès que $<\rho^{\si{\bar S}}_{\si{\bar T}},\varphi>$ l'est; autrement dit dès que la $\varphi$-composante de $\,\C^{\si{\bar S}}_{\si{\bar T}}(K_\%)$ est un $\Lambda_\varphi$-module de torsion; auquel cas on a l'estimation asymptotique stricte:\smallskip

\centerline{$x^{\si{\bar S}}_{\si{\bar T}}(n)_\varphi \ \approx \ <\mu^{\si{\bar S}}_{\si{\bar T}}, \varphi>\ell^n \ + \ <\lambda^{\si{\bar S}}_{\si{\bar T}}, \varphi>n $.}
\end{itemize}
\end{Prop}

\Preuve C'est l'application directe du Scolie 6 de \cite{J52} aux $\varphi$-composantes du module $\,\C^{\si{\bar S}}_{\si{\bar T}}(K_\%)$.\medskip

Supposons désormais $K/F$ à conjugaison complexe, $K$ contenant $\mmu^\ph_\ell$ et prenons pour $F_\%$ la $\Zl$-extension cyclotomique de $F$.
Notons $\omega$ le caractère de l'action de $\Delta$ sur $\mmu_{\ell^\%}$ et $\chi \mapsto \chi^*=\omega\chi^{-{\si{1}}}$ l'involution du miroir.
Convenons de dire qu'un caractère irréductible est réel lorsqu'il prend une valeur positive sur la conjugaison complexe; qu'il est imaginaire sinon. Écrivons $\chi=\chi^{\si{\oplus}}+\chi^{\si{\ominus}}$ la décomposition d'un caractère $\chi$ en ses composantes réelle et imaginaire.
Rappelons enfin  que le caractère d'un $\Zl[\Delta]$ module noethérien $M$ est, par convention, le caractère du tensorisé $\Ql\otimes_\Zl M$.

\begin{Th}[\bf Théorème principal]\label{TP}
Soient $\ell$ un premier impair, $F$ totalement réel, $K/F$ une extension abélienne à conjugaison complexe contenant $\mmu_\ell^\ph$, de  groupe de Galois $\Delta$ d'ordre étranger à $\ell$, et $K_\%=KF_\%$ la $\Zl$-extension cyclotomique de $K$. Soient enfin $\bar S$ et $\bar T$ deux ensembles finis disjoints de places finies de $F$ dont la réunion contient l'ensemble $L$ des places au-dessus de $\ell$, puis $S=\bar S\setminus L$ et $T=\bar T\setminus L$ leurs parties modérées et $\chi_{\si{S}}=\sum_{\p\in S}\ell^{n_{\si{\p}}}\chi_{\p}$ la somme des induits
des caractères unités des sous-groupes de décomposition respectifs $\Delta_\p$ des $\p$ de $S$ dans $K/F$ comptés avec une multiplicité égale à leur indice de décomposition $\ell^{n_{\si{\p}}}$ dans $K_\%/K$.\smallskip

Dans le cas spécial  $(\bar S,\bar T)=(LS,\emptyset)$, le défaut $\kappa_{\si{\,\emptyset}}^{\si{LS}}=-1$ est l'opposé du caractère unité. Il suit:\smallskip

\centerline{$\lambda^{\si{LS}}_{\si{\,\emptyset}}=(\lambda_{\si{L}}^{\si{\emptyset}}+\chi_{\si{S}}-1)^*$.}\smallskip

Hors le cas spécial  $(\bar S,\bar T)=(LS,\emptyset)$, le  défaut $\kappa_{\si{\bar T}}^{\si{\bar S}}$ est imaginaire: $\kappa_{\si{\bar T}}^{\si{\bar S\,\oplus}}=0$ et $\kappa_{\si{\bar T}}^{\si{\bar S}}=\kappa_{\si{\bar T}}^{\si{\bar S\,\ominus}}\geqslant 0$ .\smallskip

$(i)$  Si $\bar T$ contient  $L$, on a de plus: $\kappa_{\si{LT}}^{\si{S\,\ominus}}=0$ et le défaut $\kappa_{\si{LT}}^{\si{S}}$ est nul. Il suit alors:\smallskip

\centerline{$\lambda_{\si{LT}}^{\si{S\,\ominus}}=\lambda_{\si{L}}^{\si{S\,\ominus}}=\lambda_{\si{L}}^{\si{\emptyset\,\ominus}}+(\chi_{\si{S}}^{\si{\,\oplus}}-1)^*$,}\smallskip

$(ii)$  Si $\bar S$ contient $L$, le défaut $\kappa_{\si{T}}^{\si{LS\,\ominus}}$ est le caractère de l'image semi-locale $p_{\si{L}}^\ph(\E_{\si{T}}^{\si{\ominus}})$ du $\Zl$-module construit sur les $T$-unités imaginaires du corps $K_\%$.
Sous la conjecture de Leopoldt dans $K_\%$ (i.e. dans tous les $K_n$) et pour $\bar T$ contenant $ L$, il vient ainsi (en échangeant $\bar S$ et $\bar T$):\smallskip

\centerline{$\lambda_{\si{LT}}^{\si{S\,\oplus}}=\lambda_{\si{L}}^{\si{S\,\oplus}}=\lambda_{\si{L}}^{\si{\emptyset\,\oplus}}+(\chi_{\si{S}}^{\si{L\ominus}})^*$,}\smallskip

\noindent où $\chi_{\si{S}}^{\si{L\ominus}}$ désigne le caractère du $\Zl[\Delta]$-module des $S$-unités imaginaires $L$-infinitésimales de $K_\%$.
\smallskip

$(iii)$  Pour $F=\QQ$ enfin, et $\bar S$ ne contenant pas $\ell$, le caractère de défaut $\kappa_{\si{S}}^{\si{L}}$ est donné par:\smallskip

\centerline{$\kappa_{\si{S}}^{\si{L}} = (\chi_{\si{S}}^\ph\wedge \ell^{\,\max_{\si{p\in S}}\{n_{\si{p}}\}} \chi_\reg)^{\si{\ominus}}= \sum^{\si{\ominus}}_{\varphi}\,[\, \ell^{\,\max_{\si{p\in S_{\varphi}}}\{n_{\si{p}}\}}]\varphi$,}\smallskip

\noindent où $\varphi$ décrit les caractères irréductibles imaginaires de $\Delta$ et $S_{\varphi}=\{p\in S\;|\; {\varphi}(\Delta_p)=1\}$. Il suit:\smallskip

\centerline{$\lambda_{\si{LT}}^{\si{S\,\oplus}}=\lambda_{\si{L}}^{\si{S\,\oplus}}=\lambda_{\si{L}}^{\si{\emptyset\,\oplus}}+ \sum^{\si{\ominus}}_{\varphi}\,[\,( \sum_{p \in S_{\varphi} }\ell^{n_{\si{p}}} ) -  \ell^{\,\max_{\si{p\in S_{\varphi}}}\{n_{\si{p}}\}}]\,\varphi^*$.}\smallskip

En particulier, $\lambda_{\si{L}}^{\si{S\,\oplus}}$ et $\lambda_{\si{L}}^{\si{\,\emptyset\,\oplus}} = \lambda_{\si{\emptyset}}^{\si{\emptyset\,\oplus}}$ ont même $\varphi^*$-composante dès que $S_{\varphi}$ a au plus 1 élément.
\end{Th}

La dernière assertion recoupe et prolonge les résultats de Itoh, Mizusawa et Ozaki \cite{IMO, Ito} sur le $\Zl$-rang des modules d'Iwasawa modérément ramifiés.\smallskip

Une conséquence {\em a priori} surprenante est que, même dans le cas le plus simple où $F$ est le corps des rationnels, le caractère de défaut $\kappa_{\si{\bar T}}^{\si{\bar S\,\ominus}}$ s'avère ainsi généralement non-trivial.

\newpage
\section{Conséquence des identités du miroir de Gras}

Nous nous plaçons désormais sous les hypothèses générales du Théorème principal \ref{TP}: $\ell$ est impair, $K/F$ est à conjugaison complexe, $K$ contient les racines $\ell$-ièmes de l'unité, $F_\%$ est la $\Zl$-extension cyclotomique de $F$ et $S \sqcup T$ contient l'ensemble $L$ des places au-dessus de $\ell$.\par

Dans ce contexte les identités du miroir de Gras (cf. \cite{Gra1,J43}) mettent en reflet les sous-groupes de $\ell^n$-torsion respectifs des pro-$\ell$-groupes de $S$-classes $T$-infinitésimales et de $T$-classes $S$-infinitésimales attachés aux étages finis $K_n$ de la tour cyclotomique.\smallskip

Précisons quelques notations: pour chaque place $\p_\%$ de $F_\%$, désignons par $\Delta_{\p_\%}$ le sous-groupe de décomposition de $\p_\%$ dans $\Delta=\Gal(K_\%/F_\%)$ (qui ne dépend que de la place $\p$ de $F$ au-dessous de $\p_\%$); notons $\chi^\ph_{\p_\%}$ l'induit à $\Delta$ du caractère unité de $\Delta_{\p_\%}$; posons enfin:\smallskip

\centerline{$\chi^\ph_{\si{S}}=\sum_{\p_\%\in S_\%}\chi^\ph_{\p_\%} \qquad\&\qquad \chi^\ph_{\si{T}}=\sum_{\p_\%\in T_\%}\chi^\ph_{\p_\%}$,}\smallskip

\noindent où la somme porte sur les places de $F_\%$ au-dessus de $S$ et $T$ respectivement.\smallskip

Introduisons le caractère de Teichmüller $\omega$, défini ici comme le caractère de l'action de $\Delta$ sur le module de Tate $\Tl=\varprojlim \mmu_{\ell^n}$; notons $\chi^{\si{-1}}$ le contragrédient d'un caractère $\chi$, donné par $\sigma\mapsto\chi(\sigma^{\si{-1}})$; et écrivons

\centerline{$\chi\mapsto\chi^*=\omega\chi^{\si{-1}}$}\smallskip

\noindent l'involution du miroir. Cela étant, il vient:

\begin{Th}[\bf Théorème de réflexion]
Si le nombre premier $\ell$ est impair, si $K$ contient le groupe $\mmu_\ell$ des racines $\ell$-ièmes de l'unité, et si la réunion $\bar S\bar T$ contient l'ensemble $L$ des places de $F$ au-dessus de $\ell$, les caractères intervenant dans le Théorème des paramètres satisfont les identités:\smallskip

\centerline{$\lambda^{\si{\bar T}}_{\si{\bar S}}-\kappa^{\si{\bar T}}_{\si{\bar S}} +(\chi^\ph_{\si{\bar S}}-1)\,=\,\big(\lambda^{\si{\bar S}}_{\si{\bar T}}-\kappa^{\si{\bar S}}_{\si{\bar T}} +(\chi^\ph_{\si{\bar T}}-1)\big)^*$.}
\end{Th}

\Preuve C'est exactement l'assertion $(iii)$ du Théorème 2.6 de \cite{J43}, une fois corrigés les caractères structurels $\lambda$ des défauts $\kappa$, conformément au Théorème des paramètres plus haut.

\begin{Cor}\label{C1}
Pour tout couple $(S,T)$ d'ensembles disjoints de places modérées, il vient:\smallskip

\centerline{$\lambda_{\si{T}}^{\si{SL}}=\lambda^{\si{SL}}_{\si{\,\emptyset}}=\big(\lambda_{\si{L}}^{\si{\emptyset}}+(\chi^\ph_{\si{SL}}-1)\big)^*$.}\par
En particulier, il suit:\par
\centerline{$\lambda_{\si{T}}^{\si{SL}}=\lambda_{\si{T}}^{\si{L}}+(\chi^\ph_{\si{S}}-1)^*$.}
\end{Cor}

\Preuve Les places modérées (i.e. ne divisant pas $\ell$) étant sans inertie au-dessus de $K_\%$, il vient $\,\C_{\si{T}}^{\si{SL}}(K_\%)=\C_{\si{\emptyset}}^{\si{SL}}(K_\%)$, donc $\lambda_{\si{T}}^{\si{SL}}=\lambda^{\si{SL}}_{\si{\,\emptyset}}$; puis, en vertu du cas spécial $\kappa_{\si{\,\emptyset}}^{\si{SL}}=-1$:\smallskip

\centerline{$\lambda^{\si{SL}}_{\si{\emptyset}}= \lambda^{\si{SL}}_{\si{\emptyset}} -\kappa^{\si{SL}}_{\si{\emptyset}} +(\chi^\ph_{\si{\emptyset}}-1)=\big(\lambda_{\si{SL}}^{\si{\emptyset}} -\kappa_{\si{SL}}^{\si{\emptyset}} +(\chi^\ph_{\si{SL}}-1)\big)^*$.}\smallskip
Et, $\,\C^{\si{\emptyset}}_{\si{SL}}(K_\%)$ étant un $\Lambda$-module de torsion, on a: $\rho_{\si{SL}}^{\si{\emptyset}}=0$, donc $\kappa_{\si{SL}}^{\si{\emptyset}}=0$; d'où le résultat.

\begin{Cor}\label{C2}
Pour tout couple $(S,T)$ d'ensembles disjoints de places modérées, il vient de même:\smallskip

\centerline{$\lambda^{\si{S}}_{\si{TL}}=\lambda^{\si{S}}_{\si{L}}=\lambda_{\si{L}}^{\si{\emptyset}}+(\chi^\ph_{\si{S}}-1-\kappa^{\si{L}}_{\si{S}})^*$.}
\end{Cor}

\Preuve $\,\C^{\si{S}}_{\si{L}}(K_\%)$ est un $\Lambda$-module de torsion; d'où $\rho_{\si{L}}^{\si{S}}=0$ et  $\kappa_{\si{L}}^{\si{S}}=0$; puis, comme plus haut:\smallskip

\centerline{$\lambda^{\si{S}}_{\si{TL}}+(\chi^\ph_{\si{L}}-1)=\lambda^{\si{S}}_{\si{L}}+(\chi^\ph_{\si{L}}-1)=\lambda^{\si{S}}_{\si{L}}-\kappa^{\si{S}}_{\si{L}}+(\chi^\ph_{\si{L}}-1)= \big(\lambda^{\si{L}}_{\si{S}}-\kappa^{\si{L}}_{\si{S}}+(\chi^\ph_{\si{S}}-1)\big)^*$.}\smallskip

\noindent Les places modérées étant sans inertie au-dessus de $K_\%$, il vient $\,\C_{\si{S}}^{\si{L}}(K_\%)=\C_{\si{\emptyset}}^{\si{L}}(K_\%)$ donc $\lambda^{\si{L}}_{\si{S}}=\lambda^{\si{L}}_{\si{\,\emptyset}}$:\smallskip

\centerline{$\lambda^{\si{S}}_{\si{TL}}+(\chi^\ph_{\si{L}}-1)= \big(\lambda^{\si{L}}_{\si{\,\emptyset}} -\kappa^{\si{L}}_{\si{\,\emptyset}}-(\kappa^{\si{L}}_{\si{S}}-\kappa^{\si{L}}_{\si{\,\emptyset}})+(\chi^\ph_{\si{S}}-1)\big)^*$,}

\noindent c'est à dire:

\centerline{$\lambda^{\si{S}}_{\si{TL}}+(\chi^\ph_{\si{L}}-1)= \big(\lambda^{\si{L}}_{\si{\emptyset}}-\kappa^{\si{L}}_{\si{\emptyset}}-(\chi^\ph_{\si{\emptyset}}-1)-\kappa_{\si{S}}^{\si{L}}+(\chi^\ph_{\si{S}}-1)\big)^*$,}\smallskip

\noindent via $\kappa_{\si{\emptyset}}^{\si{L}}=-1$ en vertu du cas spécial; et finalement, par $\kappa^{\si{\emptyset}}_{\si{L}}=\rho^{\si{\emptyset}}_{\si{L}}=0$:\smallskip

\centerline{$\lambda^{\si{S}}_{\si{TL}}+(\chi^\ph_{\si{L}}-1)=\lambda^{\si{\emptyset}}_{\si{L}}+(\chi^\ph_{\si{L}}-1)+(\chi^\ph_{\si{S}}-1-\kappa_{\si{S}}^{\si{L}})^*$.}\smallskip

\begin{Def}\label{Def}
Nous disons que $\kappa_{\si{S}}^{\si{L}}=(\chi^\ph_{\si{S}}-1)-(\lambda_{\si{L}}^{\si{S}}-\lambda_{\si{L}}^{\si{\emptyset}})^*$
est le caractère de défaut attaché à l'ensemble de places modérées $S$ dans l'involution du miroir.
\end{Def}

Tout le problème est alors d'évaluer précisément $\kappa_{\si{S}}^{\si{L}}$ et de l'interpréter arithmétiquement.

\newpage
\section{Étude des composantes réelles et imaginaires}

Notons $\bar\tau\in\Delta$ la conjugaison complexe. Il est habituel de dire qu'un caractère $\ell$-adique est {\em réel} lorsque tous ses facteurs {\em absolument} irréductibles prennent la valeur $+1$ en $\bar\tau$; qu'il est {\em imaginaire} lorsque tous prennent la valeur $-1$; et de décomposer chaque caractère $\ell$-adique comme somme $\chi=\chi^{\si{\oplus}}+\chi^{\si{\ominus}}$ de ses composantes réelles et imaginaires.\smallskip

Par exemple, il résulte de la finitude bien connue du défaut de Leopoldt dans la $\Zl$-extension cyclotomique (cf. e.g. \cite{Was}) que le caractère $\rho^{\si{L}}_{\si{S}}=\rho^{\si{L}}_{\si{\,\emptyset}}$ est imaginaire. En particulier, hors le cas spécial, il suit: $0\leqslant\kappa^{\si{L\,\oplus}}_{\si{S}}\leqslant\ell^e\rho^{\si{L\,\oplus}}_{\si{S}}=0$, donc, en vertu du Corollaire \ref{C1}:

\begin{Prop}\label{pm}
En dehors du cas spécial $S=\emptyset$, où l'on a $\kappa_{\si{\,\emptyset}}^{\si{L}}=-1$, la composante réelle du caractère de défaut $\kappa^{\si{L}}_{\si{S}}$ est triviale pour tout ensemble fini $S$ de places modérées:\smallskip

\centerline{$\kappa^{\si{L\,\oplus}}_{\si{S}}=0$.}\smallskip

En particulier, alors que pour tout ensemble fini $T$ disjoint de $S$ de places modérées, on a:\smallskip

\centerline{$\lambda_{\si{TL}}^{\si{S\,\ominus}}=\lambda_{\si{L}}^{\si{S\,\ominus}}=\lambda_{\si{L}}^{\si{\emptyset\,\ominus}}+(\chi_{\si{S}}^{\si{\oplus}}-1)^*$,}\smallskip

\noindent la composante imaginaire $\kappa^{\si{L\,\ominus}}_{\si{S}}$ du défaut  peut être non-triviale et l'on a seulement:\smallskip

\centerline{$ 0 \leqslant  \lambda_{\si{L}}^{\si{S\,\oplus}} -  \lambda_{\si{L}}^{\si{\emptyset\,\oplus}}=  ( \chi_{\si{S}}^*-\kappa_{\si{S}}^{\si{L}\,*})^{\si{\oplus}} \leqslant  \chi_{\si{S}}^{*\,\si{\oplus}}$.}
\end{Prop}

De façon générale, le caractère structurel $\lambda^{\si{S}}_{\si{T}}$ est indépendant des places modérées intervenant dans $T$. C'est également le cas des places sauvages pour les composantes réelles sous la conjecture de Leopoldt en vertu de la généralisation suivante d'un résultat de Greenberg (\cite{Grb1}, Prop. 1):

\begin{Th}\label{GéGb}
Pour $\bar K$ totalement réel  et sous la conjecture de Leopoldt  pour $\ell$ à chaque étage fini de la $\Zl$-tour cyclotomique $\bar K_\%/ \bar K$, le sous-module du pro-$\ell$-groupe $\,\C^S(\bar K_\%)=\varprojlim \,\Cl^S(\bar K_n)$ construit sur les places au-dessus de $\ell$ est fini pour tout ensemble fini $S$ de places modérées.\smallskip

En particulier, dans le contexte de cette note, on a l'égalité: $\lambda_{\si{L}}^{\si{S\,\oplus}}=\lambda_{\si{\,\emptyset}}^{\si{S\,\oplus}}$.
\end{Th}

\Preuve La propriété annoncée étant asymptotique, ce n'est pas restreindre la généralité que de supposer (pour cette démonstration) que les places au-dessus de $\ell$ se ramifient totalement dans $\bar K_\%/\bar K$. Or, sous cette hypothèse, les classes de rayons modulo $S$ des places de $\bar K_n$ au-dessus de $\ell$ sont invariantes par le groupe de Galois $\Gamma_{\!n}=\Gal(\bar K_n/\bar K)$. Tout revient donc à vérifier que les classes ambiges de rayons modulo $S$ restent bornées lorsqu'on monte la tour $\bar K_\%/\bar K$.\smallskip

La formule des classes ambiges, écrite pour les $\ell$-groupes de classes $S$-infinitésimales (cf. \cite{J18}, Cor. II.2.35, ou l'appendice {\em infra} pour plus de détails) donne ici:
$$
|\,\Cl^{\si{\,S}}(\bar K_n)^{\Gamma_{\!n}}|= |\,\Cl^{\si{\,S}}(\bar K)|\;\frac{\prod e_\p(\si{\bar K_n/\bar K})}{[\bar K_n:\bar K]\,\big(\E^{\si{S}}(\bar K):\E^{\si{S}}(\bar K)\cap N_{\bar K{\si{n\!}}/\!\bar K}(\R_{\bar K_{\si{n}}})\big)}
$$
Dans celle-ci le produit des indices de ramification $e_\p(\bar K_n/\bar K)$ au numérateur est tout simplement l'indice normique $\big(\,\U_{\si{L}}(\bar K): N_{ K_{\si{n\!}}/\!\bar K}(\U_{\si{L}}(\bar K_{\si{n}}))\big)$, où $\,\U_{\si{L}}$ désigne le groupe des unités semi-locales attachées aux places de $L$; et $\,\E^{\si{S}}(\bar K)$ au dénominateur est le $\ell$-groupe des unités $S$-infinitésimales.\smallskip

Maintenant, sous la conjecture de Leopoldt, $\,\E^{\si{S}}(\bar K)$ s'identifie par le morphisme canonique de semi-localisation à un sous-module d'indice fini du noyau $\,\U_{\si{L}}^*(\bar K)$ de la norme $N_{\bar K/\QQ}$ dans $\,\U_{\si{L}}(\bar K)$. Nous pouvons donc, à un borné près, remplacer le quotient à droite par
$$
\frac{\big(\,\U_{\si{L}}(\bar K): N_{\bar K_{\si{n\!}}/\!\bar K}(\U_{\si{L}}(\bar K_{\si{n}}))\big)}{[\bar K_n:\bar K]\,\big(\,\U^*_{\si{L}}(\bar K): N_{\bar K_{\si{n\!}}/\!\bar K}(\U^*_{\si{L}}(\bar K_{\si{n}}))\big)}=1.
$$

\begin{Cor}\label{kappa}
En dehors du cas spécial et sous la conjecture de Leopoldt dans $K_\%$, le caractère de défaut $\kappa_{\si{S}}^{\si{L}}$ est donné par:

\centerline{$\kappa_{\si{S}}^{\si{L}} = [\chi^\ph_{\si{S}}-(\lambda^{\si{S}}_{\si{L}}-\lambda^{\si{\emptyset}}_{\si{L}})^*]^{\si{\ominus}}=  [\chi^\ph_{\si{S}}-(\lambda^{\si{S}}_{\si{\emptyset}}-\lambda^{\si{\emptyset}}_{\si{\emptyset}})^*]^{\si{\ominus}}$,}\smallskip

\noindent où $\lambda^{\si{\emptyset}}_{\si{\emptyset}}$ est le caractère structurel attaché à la limite projective des $\ell$-groupes de classes $\,\Cl(\bar K_n)$ des sous-corps totalement réels $\bar K_n$ et $\lambda^{\si{S}}_{\si{\emptyset}}$ à celle des $\ell$-groupes $ \,\Cl^{\si{\,S}}(\bar K_n)$ de classes de rayons modulo $S$.
\end{Cor}

\newpage
\section{Interprétation galoisienne du caractère de défaut}

Pour interpréter $\kappa^{\si{S}}_{\si{L}}$, appuyons-nous sur la théorie $\ell$-adique du corps de classes (cf. \cite{J31}).
Rappelons que, pour tout ensemble fini $\Sigma$ de places d'un corps de nombres $N$, le pro-$\ell$-groupe d'idèles associé à la pro-$\ell$-extension abélienne $\Sigma$-ramifiée maximale $H^{\si{\Sigma}}(N)$ de $N$ est le produit $\U^{\si{\Sigma}}\R$, avec\smallskip

\centerline{$\U^{\si{\Sigma}}=\prod_{\p\notin \Sigma}\,\U_\p$ \quad et \quad  $\R=\Zl\otimes_\ZZ N^\times$}\smallskip

\noindent où $\,\U_\p=\varprojlim U_\p/U_\p^{\ell^{\si{k}}}$ est le pro-$\ell$-groupe des unités en $\p$ et $\R$ celui des idèles principaux.
Posant:\smallskip

\centerline{$\U_{\si{\Sigma}}=\prod_{\p\in \Sigma}\,\U_\p$.}\smallskip

\noindent on obtient l'isomorphisme:
$\Gal(H^{\si{LS}}(N)/H^{\si{S}}(N)) \simeq \,\U^{\si{S}}\R/\U^{\si{LS}}\R  \simeq \,\U_{\si{L}}/(\,\U_{\si{L}}\cap\,\U^{\si{LS}}\R)
\simeq \,\U_{\si{L}}/p_{\si{L}}(\E^{\si{S}})$.\smallskip

\noindent Dans celui-ci $p_{\si{L}}$ est le morphisme canonique de semi-localisation, $\E=\Zl\otimes_\ZZ E$ est le $\ell$-adifié du groupe des unités et $\,\E^{\si{S}}=\E\cap\,\U^{\si{S}}$ son sous-groupe $S$-infinitésimal. De façon semblable, il vient:\smallskip

\centerline{$\Gal(H^{\si{LS}}(N)/H^{\si{L}}(N)) \simeq \,\U^{\si{L}}\R/\U^{\si{LS}}\R  \simeq \,\U_{\si{S}}/(\,\U_{\si{S}}\cap\,\U^{\si{LS}}\R)
\simeq \,\U_{\si{S}}/p_{\si{S}}(\E^{\si{L}}) \simeq \,\U_{\si{S}}$,}\smallskip

\noindent avec ici $\,\E^{\si{L}}=1$, dès que le corps $N$ vérifie la conjecture de Leopoldt en $\ell$ (cf. \cite{J31}, \S2.3). Et enfin:\smallskip

\centerline{$\Gal(H^{\si{LS}}(N)/H^{\si{L}}(N)H^{\si{S}}(N)) \simeq (\,\U^{\si{L}}\R\cap\,\U^{\si{S}}\R/\U^{\si{LS}}\R)  
\simeq (\,\U_{\si{S}}\cap \,\U^{\si{S}}\R)/(\,\U_{\si{S}}\cap\,\U^{\si{LS}}\R)
\simeq p_{\si{S}}(\E)$,}\smallskip

\noindent toujours sous la conjecture de Leopoldt dans $N$.\medskip

Appliquant cela aux étages finis $ \bar K_n$ de la $\Zl$-extension cyclotomique du sous-corps réel $ \bar K$ de $K$ et passant à la limite projective pour la norme, nous obtenons les isomorphismes de $\Lambda[\Delta]$-modules:

\begin{center}$\Gal (H^{\si{LS}}( \bar K_\%) /H^{\si{L}}( \bar K_\%) ) \simeq \overset{_\leftarrow}{\U}_{\si{S}}$,\qquad
$\Gal (H^{\si{LS}}( \bar K_\%) / H^{\si{L}} (\bar K_\%)H^{\si{S}}( \bar K_\%)) \simeq p_{\si{S}}(\overset{_\leftarrow}{\E})$,
$\Gal ( H^{\si{LS}}(\bar K_\%) /H^{\si{S}}(\bar K_\%) ) \simeq \overset{_\leftarrow}{\U}_{\si{L}}/p_{\si{L}}(\overset{_\leftarrow}{\E}{}^{\si{S}})$,\end{center}\smallskip

\noindent où $\overset{_\leftarrow}{\E} = \varprojlim \,\E( \bar K_n)$ est la limite projective des groupes d'unités; $\overset{_\leftarrow}{\U}_{\si{L}} = \varprojlim \,\U_{\si{L}}( \bar K_n)$ celle des groupes d'unités locales attachées aux places au-dessus de $\ell$
; et $\,\overset{_\leftarrow}{\U}_{\si{S}} = \varprojlim \,\U_{\si{S}}( \bar K_n) \simeq \prod_{\p_\%\in S}\,\varprojlim\, \mu_{\p_n}$ est un $\Zl[\Delta]$-module projectif de caractère $\omega\chi_{\si{S}}^\ph =\omega\chi_{\si{S}}^{\si{-1}} =\chi_{\si{S}}^*$.\smallskip

Écrivant alors $\,\C^{\si{L}}=\,\C^{\si{L}} (\bar K_\%) \simeq \,\C^{\si{L\!}}(K_\%)^{\si{\oplus}}$,
$\,\C=\,\C^{\si{\emptyset}}(\bar K_\%) \simeq \,\C^{\si{\emptyset\!}}(K_\%)^{\si{\oplus}}$,
$\,\C^{\si{S}}=\,\C^{\si{S\!}}(\bar K_\% ) \simeq \,\C^{\si{S\!}}(K_\%)^{\si{\oplus}}$, nous obtenons le diagramme galoisien:

\begin{center}
\setlength{\unitlength}{2pt}
\begin{picture}(104,108)(0,0)

\thinlines 

\put(45,102){$H^{\si{LS}}(\bar K_\%)$}
\put(37,82){$H^{\si{L}}(\bar K_\%)H^{\si{S}}(\bar K_\%)$}
\put(5,52){$H^{\si{L}}(\bar K_\%)$}
\put(87,52){$H^{\si{S}}(\bar K_\%)$}
\put(45,22){$H^{\si{\emptyset}}(\bar K_\%)$}
\put(51,0){$ \bar K_\%$}

\bezier{60}(43,102)(18,92)(13,62)
\bezier{60}(66,102)(87,92)(93,62)
\bezier{60}(13,47)(18,22)(46,4)
\bezier{60}(93,47)(88,22)(60,4)

\put(53,99){\line(0,-1){10}}
\put(53,19){\line(0,-1){10}}
\put(49,79){\line(-3,-2){32}}
\put(58,79){\line(3,-2){33}}
\put(50,28){\line(-3,2){33}}
\put(58,28){\line(3,2){33}}

\put(54,92){$p_{\si{S}}(\overset{_\leftarrow}{\E})$}
\put(12,82){$\overset{_\leftarrow}{\U}_{\si{S}}$}
\put(88,82){$\overset{_\leftarrow}{\U}_{\si{L}}/p_{\si{L}}(\overset{_\leftarrow}{\E}{}^{\si{S}})$}
\put(16,22){$\C^{\si{L}}$}
\put(86,22){$\C^{\si{S}}$}
\put(55,12){$\C$}

\end{picture}
\end{center}

\noindent où chacun des six groupes de Galois représentés est un $\Lambda|\Delta]$-module noethérien. De l'isomorphisme $\Gal (H^{\si{S}} ( \bar K_\%)/ H^{\si{\emptyset}}(\bar K_\%)) \simeq \overset{_\leftarrow}{\U}_{\si{S}}/p_{\si{S}}(\overset{_\leftarrow}{\E})$, nous tirons donc en vertu du Corollaire \ref{kappa}:

\begin{Prop}\label{kappa*}
Le reflet $\kappa_{\si{S}}^{\si{L}*}$ du caractère de défaut $\kappa_{\si{S}}^{\si{L}}$  est le caractère du $\Zl[\Delta]$-module projectif $\Gal(H^{\si{LS}}(\bar K_\%)/H^{\si{L}}(\bar K_\%)H^{\si{S}}(\bar K_\%))\simeq p_{\si{S}}(\overset{_\leftarrow}{\E}) \subset \overset{_\leftarrow}{\U}_{\si{S}}$.
\end{Prop}

\newpage
\section{Diagrammes pseudo-exacts}

\begin{Def}
Nous disons que deux sous-corps $A$ et $B$ de $\bar\QQ$ sont pseudo-équivalents, ce que nous écrivons $A \asymp B$, lorsqu'ils sont de degré fini sur un même sous-corps; autrement dit lorsque leur compositum $AB$ est de degré fini sur leur intersection: $[AB:A\cap B]<\infty$.\par
\end{Def}

\Remarque La pseudo-équivalence $A \asymp B$  définit une partition de l'ensemble des sous-corps de $\bar\QQ$:\par

-- La classe de $\QQ$ est l'ensemble des extensions finies de $\QQ$, i.e. l'ensemble des corps de nombres.\par

-- La classe de $\QQ_\%$ est formée des extensions finies de $\QQ_\%$, i.e.  des extensions $N_\%=N\QQ_\%$ où $N$ est un corps de nombres: ce sont exactement les corps surcirculaires au sens de \cite{J25}.\par

\Exemple Le Théorème \ref{GéGb} affirme sous la conjecture de Leopoldt que, pour tout ensemble fini $S$ de places modérées de $\bar K$ totalement réel, la pro-$\ell$-extension abélienne $S$-ramifiée maximale $H^{\si{S}}_{\si{\emptyset}}(\bar K_\%)$ de $\bar K_\%$ est pseudo-équivalente à sa sous-extension $L$-décomposée: {$H^{\si{S}}_{\si{\emptyset}}(\bar K_\%)\asymp H^{\si{S}}_{\si{L}}(\bar K_\%)$.

\begin{Def}
Nous disons enfin qu'un quadruplet $(A,B,C,D)$ de sous-corps de $\bar\QQ$ est pseudo-exact lorsqu'on a simultanément: $A \asymp BC$ et $D\asymp B \cap C$.

\end{Def}

Revenons maintenant au contexte qui nous intéresse: $\ell$ est un nombre premier impair; $K/F$ une extension abélienne totalement imaginaire à conjugaison complexe d'un corps totalement réel de degré $[K:F]$ étranger à $\ell$ qui contient le $\ell$-groupe $\mmu_\ell$ des racines $\ell$-ièmes de l'unité. Et plaçons-nous du point de vue de la Théorie de Kummer  (cf. \cite{J18}, I.2 ou \cite{J23}).\smallskip

Par le miroir, le radical $\bar\r_\%=\Rad(\bar K^\ab_\%K_\%/K_\%)$ de la pro-$\ell$-extension abélienne maximale $\bar K^\ab_\%$ du sous-corps réel $\bar K_\%$ de $K_\%$ est la composante imaginaire du radical  $\r_\%=(\Ql/\Zl)\otimes_\ZZ K_\%^\times$.\smallskip

Fixons un ensemble fini non vide $S$ de places modérées; et notons $\,\E_{\si{S}}^{\si{\ominus}}$ le $\Zl$-module construit sur les $S$-unités imaginaires de $K_\%$; puis $\,\E_{\si{S}}^{\si{L\,\ominus}}$ le sous-module {\em pseudo-infinitésimal} formé des $\varepsilon$ dont l'image semi-locale $p_{\si{L}}(\varepsilon)$ dans $\,\U_{\si{L}}=\prod_{\l_\%\in {\si{L}}}\U_{\l_\%}$ tombe dans  $\mu_{\si{L}}^\ph=\prod_{\l_\%\in {\si{L}}}\mu^\ph_{\l_\%}$. Observons que $\,\E_{\si{S}}^{\si{\ominus}}$ et son image  $p^\ph_{\si{L}}(\E_{\si{S}}^{\si{\ominus}})$, se lisent à un étage fini $K_{n_{\si{0}}}$ de la tour cyclotomique (cf. section 6, {\em Remarque}).

\begin{Lem}\label{LC}
Notons $H^{\si{\Sigma}}_{\si{\mathrm{T}}}(\bar K_\%)$ la $\ell$-extension abélienne $\Sigma$-ramifiée $\mathrm{T}$-décomposée maximale de $\bar K_\%$.\smallskip
\begin{itemize}
\item[(i)] Le quadruplet $(H^{\si{LS}}_{\si{\emptyset}}(\bar K_\%)K_\%,  H^{\si{L}}_{\si{\emptyset}}(\bar K_\%)K_\%, K_\%[\!\sqrt[\ell^\%]{\E_{\si{S}}^{\si{\,\ominus}}}\,], K_\%)$ est pseudo-exact.\smallskip

\item[(ii)] Le quadruplet $(H^{\si{S}}_{\si{L}}(\bar K_\%)K_\%, H^{\si{\emptyset}}_{\si{L}}(\bar K_\%)K_\%, K_\%[\!\sqrt[\ell^\%]{\E_{\si{S}}^{\si{L\,\ominus}}}\,], K_\%)$ est pseudo-exact.
\end{itemize}
\end{Lem}

\Preuve Comme $\,\E_{\si{S}}^{\si{\ominus}}$ est le produit direct de $\mmu_{\ell^\%}$ et d'un $\Zl[\Delta]$-module projectif de caractère $\chi^{\si{\ominus}}_{\si{S}}$, le groupe de Galois  $\Gal(K_\%[\!\sqrt[\ell^\%]{\E_{\si{S}}^{\si{\,\ominus}}}\,]/K_\%) \simeq \Hom ((\Ql/\Zl)\otimes_\Zl\E_{\si{L}}^{\si{\ominus}},\mmu_{\ell^\%})$ est lui-même un $\Zl[\Delta]$-module projectif de caractère $\chi^{\si{\ominus}\,*}_{\si{S}}$. Si donc la sous-extension maximale de $K_\%[\!\sqrt[\ell^\%]{\E_{\si{S}}^{\si{\,\ominus}}}]$ qui est non-ramifiée aux places de $S$ était de degré infini sur $K_\%$, elle contiendrait une $\Zl$-extension de la forme $K_\%[\!\sqrt[\ell^\%]{\varepsilon}]$ pour un certain $\varepsilon$ de $\,\E_{\si{S}}^{\si{\ominus}}$.
Or, la factorisation de l'image $(\varepsilon)$ de $\varepsilon$ dans le $\ell$-adifié $\D_{K_\%}=\Zl\otimes_\ZZ D_{K_\%}$ du groupe des diviseurs de $K_\%$ (en fait dans celui de $K_{n_{\si{0}}}$) ferait alors intervenir l'une au moins  $\p_{n_{\si{0}}}$ des places de $S$, inerte dans $K_\%/K_{n_{\si{0}}}$, laquelle serait ainsi ramifiée dans $K_\%[\!\sqrt[\ell^{\si{n}}]{\varepsilon}]/K_\%$ pour tout $n$ assez grand; et finalement dans  $H^{\si{L}}_{\si{\emptyset}}(\bar K_\%)K_\% \cap K_\%[\!\sqrt[\ell^\%]{\E_{\si{S}}^{\si{\,\ominus}}}\,]/K_\%$, ce qui est absurde.
D'où:\smallskip

\centerline{$ H^{\si{L}}_{\si{\emptyset}}(\bar K_\%)K_\% \cap K_\%[\!\sqrt[\ell^\%]{\E_{\si{S}}^{\si{\,\ominus}}}\,] \asymp K_\%$.}\smallskip

\noindent Ainsi $\Gal(H^{\si{L}}_{\si{\emptyset}}(\bar K_\%)  K_\%[\!\sqrt[\ell^\%]{\E_{\si{S}}^{\si{\,\ominus}}}\,]/ K_\%) \sim \Gal(H^{\si{L}}_{\si{\emptyset}}(\bar K_\%) K_\%/  K_\%)\times\Gal(K_\%[\!\sqrt[\ell^\%]{\E_{\si{S}}^{\si{\,\ominus}}}\,]/ K_\%)$ est un $\Lambda[\Delta]$-module de même caractère structurel $\lambda^{\si{L\,\oplus}}_{\si{\,\emptyset}}+\chi^{\si{\ominus\,}*}_{\si{S}}$ que $\Gal(H^{\si{LS}}_{\si{\emptyset}}(\bar K_\%) K_\%/ K_\%)$ par le Corollaire  \ref{C1}.
Comme il a même caractère structurel $\mu^{\si{L\,\oplus}}_{\si{\,\emptyset}}=\mu^{\si{LS\,\oplus}}_{\si{\,\emptyset}}$ d'après \cite{J43}, Th. 2.8, il lui est pseudo-isomorphe. De l'inclusion immédiate $H^{\si{L}}_{\si{\emptyset}}(\bar K_\%) K_\%[\!\sqrt[\ell^\%]{\E_{\si{S}}^{\si{\,\ominus}}}\,] \subset H^{\si{LS}}_{\si{\emptyset}}(\bar K_\%)K_\%$, on déduit donc l'équivalence:\smallskip

\centerline{$H^{\si{L}}_{\si{\emptyset}}(\bar K_\%) K_\%[\!\sqrt[\ell^\%]{\E_{\si{S}}^{\si{\,\ominus}}}\,] \asymp H^{\si{LS}}_{\si{\emptyset}}(\bar K_\%)K_\%$.}\smallskip

\noindent Des inclusions $ H^{\si{\emptyset}}_{\si{L}}(\bar K_\%)K_\% \subset  H^{\si{L}}_{\si{\emptyset}}(\bar K_\%)K_\%$ et $K_\%[\!\sqrt[\ell^\%]{\E_{\si{S}}^{\si{L\,\ominus}}}\,] \subset K_\%[\!\sqrt[\ell^\%]{\E_{\si{S}}^{\si{\,\ominus}}}\,]$, on tire par ailleurs:\smallskip

\centerline{$ H^{\si{\emptyset}}_{\si{L}}(\bar K_\%)K_\% \cap K_\%[\!\sqrt[\ell^\%]{\E_{\si{S}}^{\si{L\,\ominus}}}\,] \asymp K_\%$.}\smallskip

\noindent Enfin, puisque $\Gal(H^{\si{S}}_{\si{L}}(\bar K_\%)/\bar K_\%)$ et $\Gal(H^{\si{\emptyset}}_{\si{L}}(\bar K_\%)/\bar K_\%)$ ont même  paramètre structurel $\mu^{\si{\emptyset\,\oplus}}_{\si{\,\emptyset}}$ (toujours par \cite{J43}), pour établir la pseudo-équivalence $ H^{\si{S}}_{\si{L}}(\bar K_\%)K_\% \asymp  H^{\si{\emptyset}}_{\si{L}}(\bar K_\%)K_\%[\!\sqrt[\ell^\%]{\E_{\si{S}}^{\si{L\,\ominus}}}\,]$
il suffit de vérifier qu'une $\Zl$-extension $K_\%[\!\sqrt[\ell^\%]{\varepsilon}]$ construite sur un élément imaginaire $\varepsilon$ de $\R_\%=\Zl\otimes_\ZZ K^\times_\%$ est $S$-ramifiée et $L$-décomposée si et seulement si $\varepsilon$ est dans $\,\E_{\si{S}}^{\si{L\,\ominus}}$. Or,  la dernière condition revient à imposer que $\varepsilon$ soit $L$-pseudo-infinitésimal; et la première que ce soit une $S$-unité.

\newpage
\section{Interprétation kummérienne du caractère de défaut}

Considérons le diagramme galoisien ci-après, où apparaissent les composita respectifs avec $K_\%$ des pro-$\ell$-extensions abéliennes maximales 
$LS$-ramifiée $H^{\si{LS}}_{\si{\emptyset}}(\bar K_\%)$,
$L$-ramifiée $H^{\si{L}}_{\si{\emptyset}}(\bar K_\%)$,
$S$-ramifiée et $L$-décomposée $H^{\si{S}}_{\si{L}}(\bar K_\%)$, 
non-ramifiée et $L$-décomposée $H^{\si{\emptyset}}_{\si{L}}(\bar K_\%)$
du sous-corps réel $\bar K_\%$, ainsi que celles définies kumériennement sur $K_\%$ par le $\Zl$-module des $S$-unités imaginaires $\,\E_{\si{S}}^{\si{\,\ominus}}$ ou son sous-module $L$-pseudo-infinitésimal $\,\E_{\si{S}}^{\si{L\,\ominus}}$:

\begin{displaymath}
\xymatrix{
 & & H^{\si{LS}}_{\si{\emptyset}}(\bar K_\%)K_\%  \ar@{-}[dl]_{\kappa^{\si{L}*}_{\si{S}}} \ar@{-}[dr] & &\\
 & \mkern-36mu\mkern-18mu  
 H^{\si{L}}_{\si{\emptyset}}(\bar K_\%) H^{\si{S}}_{\si{L}}(\bar K_\%)K_\% \mkern-36mu\mkern-18mu \ar@{-}[dl] \ar@{-}[dr] &  &  
 \mkern-36mu\mkern-18mu 
 H^{\si{S}}_{\si{L}}(\bar K_\%)K_\%[\!\sqrt[\ell^\%]{\E_{\si{S}}^{\si{\,\ominus}}}\,] \mkern-18mu \ar@{-}[dl] \ar@{-}[dr]  & \\
 H^{\si{L}}_{\si{\emptyset}}(\bar K_\%)K_\% \ar@{-}[dr] & & H^{\si{S}}_{\si{L}}(\bar K_\%)K_\%   \ar@{-}[dl] \ar@{-}[dr] & &  K_\%[\!\sqrt[\ell^\%]{\E_{\si{S}}^{\si{\,\ominus}}}\,] \ar@{-}[dl]^{\;(\chi_{\si{S}}^{\si{\ominus}}-\chi^{\si{L\ominus}}_{\si{S}})^*} \\
 &  H^{\si{\emptyset}}_{\si{L}}(\bar K_\%)K_\% \ar@{-}[dr] & &  K_\%[\!\sqrt[\ell^\%]{\E_{\si{S}}^{\si{L\,\ominus}}}\,] \ar@{-}[dl]^{(\chi^{\si{L\ominus}}_{\si{S}})^*} & \\
& &  K_\% & &
}
\end{displaymath}\bigskip

Il résulte du Lemme \ref{LC} que les quadruplets des sommets de chacun des cinq losanges représentés sont pseudo-exacts. En particulier, nous avons donc un pseudo-isomorphisme de $\Lambda[\Delta]$-modules:\smallskip

\centerline{$\Gal(H^{\si{LS}}_{\si{\emptyset}}(\bar K_\%)K_\%/H^{\si{L}}_{\si{\emptyset}}(\bar K_\%)H^{\si{S}}_{\si{L}}(\bar K_\%)K_\%) \sim \Gal(K_\%[\!\sqrt[\ell^\%]{\E_{\si{S}}^{\si{\,\ominus}}}\,]/K_\%[\!\sqrt[\ell^\%]{\E_{\si{S}}^{\si{L\,\ominus}}}\,])$.}\smallskip

\noindent Or, dans le groupe de Galois à gauche, raisonnant toujours à pseudo-isomorphisme près, nous pouvons remplacer $H^{\si{S}}_{\si{L}}(\bar K_\%)$ par $H^{\si{S}}_{\si{\emptyset}}(\bar K_\%)$ sous la conjecture de Leopoldt dans $\bar K_\%$ en vertu du Théorème \ref{GéGb}, comme expliqué dans la section précédente. Il vient donc finalement:\smallskip

\centerline{$\Gal(H^{\si{LS}}_{\si{\emptyset}}(\bar K_\%)/H^{\si{L}}_{\si{\emptyset}}(\bar K_\%)H^{\si{S}}_{\si{\emptyset}}(\bar K_\%)) \sim \Gal(K_\%[\!\sqrt[\ell^\%]{\E_{\si{S}}^{\si{\,\ominus}}}\,]/K_\%[\!\sqrt[\ell^\%]{\E_{\si{S}}^{\si{L\,\ominus}}}\,])$.}\smallskip

\noindent Maintenant, le groupe de Galois à gauche est un $\Zl[\Delta]$-module projectif de caractère $\kappa_{\si{L}}^{\si{S}*}$, en vertu de la Proposition \ref{kappa*}. Quant au caractère du groupe de Galois à droite, c'est tout simplement le reflet $(\chi_{\si{S}}^{\si{\ominus}}-\chi^{\si{L\ominus}}_{\si{S}})^*$ du caractère du $\Zl[\Delta]$-module quotient $\,\E_{\si{S}}^{\si{\,\ominus}}/\E_{\si{S}}^{\si{L\,\ominus}}$, lequel s'identifie canoniquement à l'image semi-locale $p^\ph_{\si{L}}(\E_{\si{S}}^{\si{\ominus}})$ modulo torsion du $\Zl[\Delta]$-module $\E_{\si{S}}^{\si{\ominus}}$ des $S$-unités imaginaires de $K_\%$ dans le groupe $\,\U_{\si{L}}/\mu_{\si{L}}^\ph=\prod_{\l_\%\in{\si{L}}}\,\U^\ph_{\l_\%}/\mu^\ph_{\l_\%}$. Ainsi:

\begin{Th}\label{ThKappa}
Sous la conjecture de Leopoldt dans $K_\%$, le caractère de défaut $\kappa^{\si{L}}_{\si{S}}$  est le caractère $\chi_{\si{S}}^{\si{\ominus}}-\chi^{\si{L\ominus}}_{\si{S}}$ du $\Zl[\Delta]$-module $\,\E_{\si{S}}^{\si{\ominus}}/\E_{\si{S}}^{\si{L\,\ominus}}$; i.e. de l'image semi-locale $p^\ph_{\si{L}}(\E_{\si{S}}^{\si{\ominus}})$ (modulo $\Zl$-torsion) du $\Zl[\Delta]$-module $\E_{\si{S}}^{\si{\ominus}}$ des $S$-unités imaginaires de $K_\%$.
\end{Th}

\Remarque Un intérêt essentiel de ce résultat est que, comme observé plus haut, le pro-$\ell$-groupe $\E_{\si{S}}^{\si{\ominus}}$ des $S$-unités imaginaires de $K_\%$, et donc son image  $p^\ph_{\si{L}}(\E_{\si{S}}^{\si{\ominus}})$,  se lisent à un étage fini $K_{n_{\si{0}}}$ de la tour cyclotomique; en fait dès que celles des places au-dessus de $S$ qui sont décomposées par la conjugaison complexe ne se décomposent pas dans $K_\%/K_{n_{\si{0}}}$.

\begin{Cor}
Sous la conjecture de Leopoldt dans $\bar K_\%$, pour tout couple $(S,T)$ d'ensembles finis disjoints de places modérées, on a les identités entre caractères structurels:\smallskip

\centerline{$\lambda_{\si{TL}}^{\si{S\,\ominus}}=\lambda_{\si{L}}^{\si{S\,\ominus}}=\lambda_{\si{L}}^{\si{\emptyset\,\ominus}}+(\chi_{\si{S}}^{\si{\oplus}}-1)^*$\qquad et \qquad $\lambda_{\si{TL}}^{\si{S\,\oplus}}=\lambda_{\si{L}}^{\si{S\,\oplus}}=\lambda_{\si{L}}^{\si{\emptyset\,\oplus}}+(\chi_{\si{S}}^{\si{L\ominus}})^*$.}\smallskip

\noindent où $\chi_{\si{S}}^{\si{L\ominus}}$ est le caractère  du $\Zl[\Delta]$-module $\E_{\si{S}}^{\si{L\ominus}}$ des $S$-unités imaginaires $L$-infinitésimales de $K_\%$.
\end{Cor}

\newpage
\section{Application aux corps totalement réels}

Partons d'un corps totalement réel arbitraire $F$, posons $K=F[\zeta_\ell]$ et notons $\bar K=F[\zeta_\ell+\bar\zeta_\ell]$ son sous-corps réel.
Nous supposons par commodité dans ce qui suit que $F$ est galoisien sur $\QQ$, mais ce n'est pas vraiment une restriction: si tel n'est pas le cas, il suffit de remplacer $F$ par sa clôture galoisienne $F'$, puis de redescendre les résultats sur $F$ à l'aide de la norme $N_{\si{F'/F}}$.\smallskip

Étant donné un ensemble fini $S$ de places modérées, stable par $\Gal(F/\QQ)$, autrement dit un ensemble fini de nombres premiers $p\ne\ell$, désignons par $K_{n_{\si{0}}}$ le plus petit étage de la $\Zl$-tour $K_\%$ au-dessus duquel ceux-ci ne se décomposent plus et par $G_{n_{\si{0}}}$ le groupe $\Gal (K_{n_{\si{0}}}/\QQ)$. Le $\Zl$-module $\,\E_{\si{S}}^{\si{\,\ominus}}$ des $S$-unités imaginaires de $K_\%$, regardé modulo le sous-groupe de torsion, est alors un $\Zl[G_{n_{\si{0}}}]$-module projectif de caractère $\psi_{\si{S}}^{\si{\ominus}} =\sum_{p\in S} \psi_p^{\si{\ominus}}$, où $\psi_p^\ph =\Ind^{G_{n_{\si{0}}}}_{D_p}1_{D_p}$ est l'induit à $G_{n_{\si{0}}}$ du caractère unité du sous-groupe de décomposition $D_p$ de $p$ dans $K_{n_{\si{0}}}/\QQ$.\smallskip

Cela étant, sous la conjecture générale d'indépendance $\ell$-adique énoncée dans \cite{J10}, établie pour les corps abéliens et appliquée dans $K_{n_{\si{0}}}$, le caractère de l'image semi-locale $p^\ph_{\si{L}}(\E_{\si{S}}^{\si{\ominus}})$ est le plus grand caractère $\psi_{\si{S}}^{\si{L\,\ominus}}=\psi_{\si{S}}^{\si{\ominus}}\wedge \psi^\reg_{G_{n_{\si{0}}}}$ de $G_{n_{\si{0}}}$ contenu dans $\psi_{\si{S}}^{\si{\ominus}}$ comme dans le caractère régulier $\psi^\reg_{G_{n_{\si{0}}}}$.\smallskip

La conjecture $\ell$-adique entraînant celle de Leopoldt, par restriction à $G = \Gal(K/\QQ)$ il suit:

\begin{Th}\label{Calcul}
Soit $F$ un corps réel absolument galoisien, $\ell$ un nombre premier impair, $K=F[\zeta_\ell]$ et $\bar K=F[\zeta_\ell+\bar\zeta_\ell]$ son sous-corps réel. Étant donné un ensemble fini arbitraire $S$ de nombres premiers $p\ne\ell$, pour chaque $p\in S$ notons $K_{n_{\si{p}}}$ le sous-corps de décomposition de $p$ dans la tour cyclotomique $K_\%/K$ et $G_p$ son sous-groupe de décomposition dans $G=\Gal(K/\QQ$).\smallskip

Sous la conjecture d'indépendance $\ell$-adique pour $K_\%$ avancée dans \cite{J10}, et donc inconditionnellement pour $F$ abélien, le caractère du $\Zl[G]$-module $p^\ph_{\si{L}}(\E_{\si{S}}^{\si{\ominus}})$ est donné par: $\sum^{\si{\ominus}}_{\varphi}\ell^{\,\max_{\si{p\in S_{\varphi}}}\{n_{\si{p}}\}}\, \varphi$.\smallskip

Dans cette formule $\varphi$ décrit les caractères $\ell$-adiques irréductibles imaginaires de $G$ et, pour $\varphi$ fixé, $p$ décrit le sous-ensemble $S_{\varphi}=\{p\in S\;|\; {\varphi}(G_p)=1\}$ des $p$ qui vérifient ${\varphi} \leqslant \chi_p$.
\end{Th}

\begin{Cor}
Prenant $F=\QQ$, on obtient inconditionnellement pour les invariants de $\bar K$:\smallskip

\centerline{$\lambda_{\si{\,\emptyset}}^{\si{S\,\oplus}}=\lambda_{\si{L}}^{\si{S\,\oplus}}=\lambda_{\si{L}}^{\si{\emptyset\,\oplus}}+ \sum^{\si{\ominus}}_{\varphi}\,[\,( \sum_{p \in S_{\varphi} }\ell^{n_{\si{p}}} ) -  \ell^{\,\max_{\si{p\in S_{\varphi}}}\{n_{\si{p}}\}}]\varphi^*$.}\smallskip

En particulier, $\lambda_{\si{L}}^{\si{S\,\oplus}}$ et $\lambda_{\si{L}}^{\si{\emptyset\,\oplus}} = \lambda_{\si{\,\emptyset}}^{\si{\emptyset\,\oplus}}$ ont même $\varphi^*$-composante dès que $S_{\varphi}$ a au plus 1 élément.
\end{Cor}

\Exemple Prenant $\varphi^*=1$, on obtient la valeur donnée par Itoh, Mizusawa et Ozaki \cite{IMO}:\smallskip

\centerline{$\lambda^{\si{S}}(\QQ_\%)=\big(\sum_{p \in S_\omega}\ell^{n_{\si{p}}}\big)-\ell^{\,\max_{\si{p\in S_\omega}} \{n_{\si{p}}\}}$,}\smallskip

\noindent où la somme et le maximum portent sur les seuls premiers $p$ de $S$ pour lesquels $\omega=1^*$ est représenté dans $\chi_p$, i.e. sur les premiers $p$ de $S$ complètement décomposés dans $\QQ[\mmu_\ell]/\QQ$.
\medskip

\Preuve Le Corollaire résulte de la Proposition \ref{pm}, compte tenu de l'égalité $\kappa^{\si{L\,\ominus}}_{\si{S}}=\chi_{\si{S}}^{\si{\ominus}}-\chi^{\si{L\ominus}}_{\si{S}}$ donnée par le Théorème \ref{ThKappa} et de l'expression de $\chi_{\si{S}}^{\si{\ominus}}-\chi^{\si{L\ominus}}_{\si{S}}$ donnée par le Théorème ci-dessus.\par
Et lorsque $S_\varphi$ est soit vide soit un singleton $\{p\}$, le terme correctif $\big(\sum_{p \in S_\omega}\ell^{n_{\si{p}}}\big) - \ell^{\,\max_{\si{p\in S_\omega}}\{n_{\si{p}}\}}$ est nul. C'est évidemment toujours le cas lorsque $S$ lui-même est un singleton.
\medskip

\Remarque
En cohérence avec l'ensemble de la note, nous avons imposé ici que le degré $[K:F]$ de l'extension abélienne considérée soit étranger à $\ell$. Mais cette restriction n'est nullement nécessaire pour invoquer le Théorème \ref{Calcul}. D'autre part, dans le diagramme de la section précédente, les divers groupes de Galois qui interviennent sont des $\Zl$-modules noethériens. Or, pour un tel module $\X$, l'invariant $\lambda$ d'Iwasawa n'est autre que la dimension du $\Ql$-espace vectoriel $V_\X=\Ql\otimes_\Zl\X$.\par

Il est donc encore loisible de décomposer $V_\X$ comme somme directe de ses composantes isotypiques $V_\X^{e_{\si{\varphi}}}$ indexées par les caractères $\ell$-adiques irréductibles de $\Delta$ et de définir le {\em caractère} $\lambda$ de $\Delta$ attaché à $\X$ par la formule $\lambda= \sum_\varphi\, \lambda_\varphi \varphi$, avec $\lambda_\varphi =\dim_{\Ql} V_\X^{e_{\si{\varphi}}}/\deg\,\varphi$, quand bien même l'hypothèse $\ell\nmid [K:F]$ serait en défaut. De ce fait, les formules obtenues pour les composantes réelles des invariants $\lambda$ sont encore valides dans ce cadre plus large.\par

On retrouve ainsi, par exemple, le fait que pour un corps abélien réel $\bar K$, on a $\lambda_{\si{T}}^{\si{S}}=\lambda_{\si{\,\emptyset}}^{\si{\emptyset}}$ lorsque $S$ est un singleton $\{p\}$ avec $p\ne\ell$, quel que soit $T$ fini ne contenant pas $p$.

\newpage
\section*{Appendice: Suite exacte des classes infinitésimales ambiges}
\addcontentsline{toc}{section}{Appendice: Identité des classes infinitésimales ambiges}

Nous reproduisons dans cet appendice pour la commodité du lecteur une preuve succincte de la formule des classes ambiges dans le cas particulier des $\ell$-classes $S$-infinitésimales qui nous intéresse ici. Nous renvoyons à \cite{J18},  II.2 pour une étude équivariante plus générale.\smallskip

Les données sont les suivantes: $\ell$ est un nombre premier impair; $N/K$ est une $\ell$-extension cyclique de groupe $\Gamma$, et $S$ est un ensemble fini de places finies de $K$ étrangères à $\ell$.\par

Pour chacun des corps ci-dessus, par exemple $N$, le $\ell$-groupe des classes $S$-infinitésimales $\,\Cl_N^{\,\si{S}}$ n'est autre que le $\ell$-adifié $\Zl\otimes_\ZZ Cl_N^{\,\m}$ du groupe des classes de rayons modulo $\m_{\si{N}}^{\si{S}}=\prod_{\p_{\si{N}}\mid S}\p^\ph_{\si{N}}$ défini comme le quotient $Cl_N^{\,\m}=D_N^{\,\m}/P_N^{\,\m}$ du groupe $D_N^{\,\m}$ des diviseurs étrangers à $\m_{\si{N}}^{\si{S}}$ par le sous-groupe $P_N^{\,\m}$ des diviseurs principaux engendrés par les $x$ de $N^\times$ qui vérifient $x\equiv 1 \;{\rm mod}^\times \m_{\si{N}}^{\si{S}}$.\smallskip

Il vient: $\Cl_N^{\,\si{S}}=\D_N^{\,\si{S}}/\P_N^{\,\si{S}}$ avec $\D_N^{\,\si{S}}=\Zl\otimes_\ZZ D_N^{\,\si{S}}$ et $\P_N^{\,\m}=\{(x)\in \D_N^{\,\m} \,|\, p_{\si{S}}(x)=1\}$, puisqu'aux places étrangères à $\ell$, les $\ell$-adifiés $\,\U_{\p^\ph_{\si{N}}}\!$ des groupes d'unités locales $U_{\p^\ph_{\si{N}}}\!$ se réduisent aux $\ell$-groupes $\mu_{\p^\ph_{\si{N}}}\!$ de racines de l'unité, de sorte que les éléments de $\R_N=\Zl\otimes_\ZZ N^\times$ construits sur les $x$ de $N^\times$ qui vérifient la congruence précédente sont précisément ceux d'image locale  triviale aux places au-dessus de $S$; i.e. les éléments du sous-groupe $S$-infinitésimal $\R_N^{\,\m}=\{x\in \R_N \,|\, p_{\si{S}}(x)=1\}$.

\begin{Th*}[Classes $S$-infinitésimales ambiges]
Soient $\ell$ un nombre premier impair, $N/K$ une $\ell$-extension cyclique de groupe $\Gamma$ et $S$ un ensemble fini de places finies de $K$ étrangères à $\ell$.\par

Alors le nombre de $\ell$-classes de $\,\Cl^{\,\si{S}}_N$ qui sont invariantes par $\Gamma$ est donné par:
$$
|\,\Cl^{\,\si{S}\,\Gamma}_N| = |\,\Cl^{\,\si{S}}_K|\;\frac{\prod_{\p_{\si{K}}\notin S} e_{\p_{\si{K\!}}}(N/K) } { [N:K]\;\big(\E^{\,\si{S}}_K:\E^{\,\si{S}}_K\cap N_{L/K}(\R^\ph_N)\big)}
$$
\noindent où $e_{\p^\ph_{\si{K\!}}}(N/K)$ est l'indice de ramification de $\p_{\si{K}}$ et $\,\E_K^{\,\si{S}}$ le groupe des unités $S$-infinitésimales.
\end{Th*}

\Preuve Elle est essentiellement identique à celle de la formule de Chevalley (\cite{Chv}, pp. 402--406).\smallskip

$(i)$ Comparaison des classes ambiges et des classes d'ambiges: on dispose d'un isomorphisme\smallskip

\centerline{$\Cl^{\,\si{S}\,\Gamma}_N / cl^{\si{S}}(\D^{\,\si{S}\,\Gamma}_N) \simeq \E^{\,\si{S}}_K\cap N_{N/K}(\R^\ph_N) / N_{L/K}(\E^{\,\si{S}}_L)$,}\smallskip

\noindent obtenu en prenant un générateur arbitraire $\sigma$ de $\Gamma$ et en envoyant la classe d'un idéal $\a$ qui vérifie $\a^{\sigma-{\si{1}}}=(\alpha)$ sur celle de l'élément $\varepsilon = N_{L/K}(\alpha)$. D'où l'identité:\smallskip

\centerline{$\big(\Cl^{\,\si{S}\,\Gamma}_N : cl^{\si{S}}(\D^{\,\si{S}\,\Gamma}_N)\big) = \frac{\big(\E^{\,\si{S}}_K: N_{N/K}(\E^{\,\si{S}}_N)\big) }{\big(\E^{\,\si{S}}_K:\;\E^{\,\si{S}}_K \cap \,N_{N/K}(\R^\ph_N) \big)}$}\smallskip

$(ii)$ Comparaison des classes d'ambiges et des classes étendues: on a l'égalité immédiate\smallskip

\centerline{$|cl^{\si{S}}(\D^{\,\si{S}\,\Gamma}_N)| =\big(\D^{\,\si{S}\,\Gamma}_N : \P^{\,\si{S}\,\Gamma}_N\big) = \frac{\big(\D^{\,\si{S}\,\Gamma}_N :\,\D^{\,\si{S}}_K\big) \; \big(\D^{\,\si{S}}_K :\,\P^{\,\si{S}}_K\big) }{\big(\P^{\,\si{S}\,\Gamma}_N :\,\P^{\,\si{S}}_K\big) }$}\smallskip

\noindent avec, au numérateur, $\big(\D^{\,\si{S}\,\Gamma}_N :\,\D^{\,\si{S}}_K\big) =  \prod_{\p_{\si{K}}\notin S} e_{\p_{\si{K\!}}}(N/K)$ et 
$ \big(\D^{\,\si{S}}_K:\P^{\,\si{S}}_K\big) =  |\Cl^{\,\si{S}}_K|$.\medskip

$(iii)$ Interprétation cohomologique du dénominateur: $\big(\P^{\,\si{S}\,\Gamma}_N :\,\P^{\,\si{S}}_K\big)=|H^1(\Gamma,\E_N^{\,\si{S}})|$\smallskip

\noindent Partant de la suite exacte courte $1 \to \E^{\si{S}}_N \to \R^{\si{S}}_N \to \P^{\si{S}}_N \to 1$, prenant ensuite la suite longue de cohomologie et la comparant à la suite de départ écrite pour $K$, on obtient la suite exacte\smallskip

\centerline{$1 \to \P^{\,\si{S}}_K \to \P^{\,\si{S}\,\Gamma}_N  \to  H^1(\Gamma,\E^{\si{S}}_N) \to  H^1(\Gamma,\R^{\si{S}}_N)$}\smallskip

\noindent Or, prenant la suite de localisation $1 \to \R_N^{\si{S}} \to \R_N \to \prod_{\p\in S} \R_{N_\p} \to 1$, puis la cohomologie, on a:\smallskip

\centerline{$1 \to \R_K^{\si{S}} \to \R_K \to \prod_{\p\in S} \R_{K_\p} \to H^1(\Gamma, \R_N^{\si{S}}) \to H^1(\Gamma, \R_N)=1$}\smallskip

\noindent et le terme de droite est trivial en vertu du Théorème 90 de Hilbert; d'où: $H^1(\Gamma, \R_N^{\si{S}})=1$.\medskip

$(iv)$ Utilisation du quotient de Herbrand $q(\Gamma,E_N)=\dfrac{|H^2(\Gamma,E_N)|}{|H^1(\Gamma,E_N)|}=\dfrac{1}{[N:K]}$ (cf. \cite{He1,He2}):\smallskip

Observant que $\,\E^{\si{S}}_N$ est d'indice fini dans $\,\E_N$ on a: $q(\Gamma,\E^{\si{S}}_N)=q(\Gamma,\E_N)=q(\Gamma,E_N)=\frac{1}{[N:K]}$.\medskip

Récapitulant le tout, on obtient le résultat annoncé.

\newpage

\section*{\normalsize\sc Addendum}

Le calcul des caractères structurels $\rho^{\si{\bar S}}_{\si{\bar T}}$ et $\mu^{\si{\bar S}}_{\si{\bar T}}$ est effectué dans \cite{J43}: le premier est purement galoisien; le second conjecturalement nul (et effectivement pour $K$ abélien).
\smallskip

L'erreur sur le module de défaut, introduite dans \cite{J18}, et reproduite dans \cite{J40} puis dans \cite{J43} a été repérée par Salle \cite{Sal} puis corrigée  dans \cite{J52} en collaboration avec Maire et Perbet. Comme expliqué dans l'introduction, le but premier de cette note est de préciser cette correction en termes de caractères en  formulant correctement des identités du miroir de Gras pour les modules d'Iwasawa.\smallskip

Les résultats présentés recoupent ceux de Itoh, Mizusawa et Ozaki \cite{IMO} ainsi que ceux de Itoh \cite{Ito} sur les modules d'Iwasawa modérément ramifiés. L'approche d'Itoh, totalement différente de la nôtre, repose sur le théorème de Kronecker-Weber et la description explicite des annulateurs pour les modules d'Iwasawa dans les tours cyclotomiques. Accessoirement elle utilise en outre les résultats de Khare et Wintenberger \cite{KW1,KW2} sur certains radicaux de Kummer.\smallskip

Je remercie enfin tout particulièrement Ch. Maire et G. Gras ainsi que le rapporteur anonyme pour leur lecture critique.

\def\refname{\normalsize{\sc  Références}}

\medskip\noindent
{\small
\begin{tabular}{l}
Institut de Mathématiques de Bordeaux \\
Université de {\sc Bordeaux} \& CNRS \\
351 cours de la libération\\
F-33405 {\sc Talence} Cedex\\
courriel : Jean-Francois.Jaulent@math.u-bordeaux.fr\\
{\footnotesize \url{https://www.math.u-bordeaux.fr/~jjaulent/}}
\end{tabular}
}

 \end{document}